\documentclass[10pt]{article}
\usepackage{amsmath,amssymb,amsthm,mathrsfs,dsfont}
\usepackage[margin=3cm]{geometry}
\usepackage{titlesec,hyperref}
\usepackage{color}
\usepackage{fancyhdr}
\usepackage{subfigure}
\usepackage[graphicx]{realboxes}
\usepackage{float}
\usepackage{empheq}
\titleformat{\subsection}
  {\it}
  {\thesubsection.\enspace}
  {1pt}
  {}

\newtheorem{definition}{Definition}[section]
\newtheorem{theorem}{Theorem}[section]
\newtheorem{lemma}{Lemma}[section]
\newtheorem{remark}{Remark}[section]
\newtheorem{proposition}{Proposition}[section]

\numberwithin{equation}{section}

\allowdisplaybreaks

\begin{document}
\date{}
\title{Local well-posedness, wave breaking and ill-posedness for the unidirectional non-local wave model}
\author{
Jian Chen$^{1}$\thanks{Corresponding author: chenj883@mail2.sysu.edu.cn (Jian Chen)},~~
Shaojie Yang$^{2}$\thanks{shaojieyang@kust.edu.cn (Shaojie Yang)},~~
Zhaoyang Yin$^{1}$\thanks{mcsyzy@mail.sysu.edu.cn (Zhaoyang Yin)}
\\[2mm]
\small $^{1}$School of Science,
\small Shenzhen Campus of Sun Yat-sen University, Shenzhen 518107, China
\\
\small $^{2}$ Department of Mathematics,  Faculty of Science, Kunming University of Science and Technology,\\
 \small Kunming,  Yunnan 650500, China
}

\date{}
\maketitle
\begin{abstract}
In this paper, we first establish local well-posedness (existence, uniqueness and continuous dependence) for the Cauchy problem associated with the unidirectional non-local wave model in the Besov spaces $B^s_{p,r}$ with $s>1+\frac{1}{p}$ and $1\leq p,r\leq\infty$, as well as in the critical Besov spaces $B^{1+\frac{1}{p}}_{p,1}$ with $1\leq p<\infty$. Next, by exploiting the fine structure of the equation and performing a delicate analysis of the commutator terms, we establish a new sufficient condition for wave breaking without using any conservation law. Finally, we prove that the model is ill-posed in $B^{1+\frac{1}{p}}_{p,r}$ for $1\leq p\leq\infty$ and $1<r\leq\infty$. Our analysis further clarifies the well-posedness theory and wave-breaking behavior of the unidirectional non-local wave model arising from the motion of a cold plasma in a magnetic field, which refines and improves previously known related results.

 \noindent\emph{Keywords}: The unidirectional non-local wave model; Local well-posedness; Wave breaking; ill-posedness
\\

\noindent\emph{Mathematics Subject Classification}: 35B44; 35R25; 35Q35
\end{abstract}
\noindent\rule{\textwidth}{0.5pt}

\section{Introduction}\label{sec1}
In 1960, Gardner and Morikawa introduced the following hyperbolic-hyperbolic-elliptic system to study the propagation of hydromagnetic waves across a magnetic field \cite{r43}:
\begin{subequations}\label{eq1}
\begin{align}
n_t+(un)_x &= 0, \label{eq1a}\\
u_t+uu_x+\frac{bb_x}{n} &= 0, \label{eq1b}\\
b-n-\left(\frac{b_x}{n}\right)_x &= 0, \label{eq1c}
\end{align}
\end{subequations}
where $n=n(t,x)$, $u=u(t,x)$, and $b=b(t,x)$ denote the ionic density, the ionic velocity, and the magnetic field of a cold plasma, respectively. This system also arises as a reduced model for collision-free two-fluid dynamics when electron inertia, charge separation, and displacement currents are neglected. Since its introduction, system \eqref{eq1} has attracted considerable attention. Gardner and Morikawa \cite{r43} formally derived a KdV approximation for small-amplitude, long-wavelength waves. Berezin and Karpman \cite{r44} studied nonstationary waves of small but finite amplitude propagating perpendicular and obliquely to the magnetic field. Taniuti and Wei \cite{r23} developed a reductive perturbation method, which Kakutani et al. \cite{r24} applied to obtain the KdV equation at leading order for all propagation angles. Pu and Li \cite{r707} rigorously justified the KdV limit. Furthermore, the local well-posedness of classical solutions with initial data $n_0(x)>0$ and $(n_0-1,u_0)\in H^2(\mathbb{R})\times H^3(\mathbb{R})$ was established in \cite{r808}. Bae, Choi, and Kwon in \cite{r909} proved the finite time $C^1$ blow up of its solutions for a class of initial data. 

Recently, by means of a multiple-scale expansion
\cite{r116,r118,r3,r23} together with a standard
unidirectional reduction procedure \cite{r12}, Alonso-Or\'{a}n, Dur\'{a}n, and Granero-Belinch\'{o}n \cite{r1} derived the following unidirectional non-local wave model from system \eqref{eq1}:
\begin{align}\label{dxfc}
h_t+\frac{3}{2}hh_x=\frac{1}{2}Lh_x+\frac{1}{2}h_x-\frac{1}{2}[L,Lh_x]h,
\end{align}
where $L=(1-\partial_x^2)^{-1}$ is the Bessel potential operator, defined by the convolution
\begin{align*}
Lf(x)=G*f(x)=\int_{\mathbb{R}}G(x-y)f(y)dy,\qquad G(x)=\frac{1}{2}\mathrm{e}^{-|x|}
\end{align*}
and $[L,f]g = L(fg) - fLg$ denotes the commutator. Eq.~\eqref{dxfc} preserves the following two conserved quantities:
\begin{align*}
\int_{\mathbb{R}}h(t,x)\,dx=\int_{\mathbb{R}}h_0(x)\,dx,
\qquad
\int_{\mathbb{R}}h^2(t,x)\,dx=\int_{\mathbb{R}}h_0^2(x)\,dx.
\end{align*}
The local well-posedness of Eq.~\eqref{dxfc} in $H^s(\mathbb{R})$ for $s>\frac{3}{2}$ was established in \cite{r1}, along with a sufficient condition for wave breaking based on a sufficiently negative initial slope, with the threshold depending on $\|h_0\|_{L^2}$ and $\|h_0\|_{L^\infty}$. Further wave breaking criteria, including a local-in-space criterion, were established in \cite{r2,r31}. More recently, it was shown that wave breaking may occur even for small initial slopes \cite{r32}. The existence of small-amplitude periodic traveling waves was proved using bifurcation theory, and other traveling waves with monotone and non-monotone profiles and varying regularity were investigated through normal-form analysis and numerical computations \cite{r303}.

The unidirectional non-local wave model \eqref{dxfc} resembles the well-known Fornberg-Whitham (FW) equation:
\begin{align}\label{eq8}
u_t+\frac{3}{2}uu_x=Lu_x,
\end{align}
which was proposed by Whitham \cite{r5} and further investigated by Fornberg and Whitham \cite{r6} as a nonlocal dispersive model to describe wave breaking phenomena, i.e., the wave profile remains bounded while its spatial derivative becomes unbounded in finite time. In recent years, the Cauchy problem for \eqref{eq8} has been extensively studied in the mathematical literature. Holmes \cite{r15} proved well-posedness of the FW equation in Sobolev spaces $H^s(\mathbb{T})$ with $s>\frac{3}{2}$. Later, Holmes and Thompson \cite{r16} established its local well-posedness in Besov spaces $B^s_{2,r}(\mathbb{R}~\text{or}~\mathbb{T})$ with $s>\frac{3}{2}$, $1<r<\infty$ or $s\geq\frac{3}{2}$, $r=1$. Guo \cite{r51} further established local well-posedness, including continuous dependence, for the FW equation in $B^s_{p,r}(\mathbb{R})$ with $s>1+\frac1p$, $1\leq p,r\leq\infty$, and in $B^{1+\frac1p}_{p,1}(\mathbb{R})$ with $1\leq p<\infty$. Qu et al. \cite{r901} subsequently established local well-posedness in the critical Besov space $B^1_{\infty,1}(\mathbb{R})$ and proved that the data-to-solution map is not uniformly continuous. Li et al. \cite{r52} established ill-posedness in $B^s_{p,r}(\mathbb{R})$ for $1<s<1+\frac1p$, $2\leq p<\infty$, and $1\leq r\leq\infty$, as well as for $s=1$, $2\leq p<\infty$, and $1\leq r\leq2$, by proving norm inflation for some specially chosen initial data. Constantin and Escher \cite{r13} first established a sufficient condition for wave breaking in terms of the minimum and maximum slopes of the initial profile, and Haziot \cite{r17} derived further slope conditions through refined estimates. Wei \cite{r19} established a sufficient condition for wave breaking and obtained estimates for the lifespan and blow-up rate, and subsequently derived new wave breaking criteria by analyzing
a Riccati-type inequality involving a time-dependent function \cite{r21}. H\"{o}rmann \cite{r18} investigated wave breaking for periodic solutions and identified classes of initial profiles leading to finite-time blow-up. Yang \cite{r22} established another sufficient condition by exploiting the structure of the equation and the $L^2$-conservation law. For further results on the FW equation, we refer the reader to
\cite{r905,r20,r902,r904,r903} and the references therein.

Although local well-posedness of Eq.~\eqref{dxfc} in $H^s(\mathbb{R})$ for $s>\frac{3}{2}$ was established in \cite{r1}, owing to the complexity of the nonlocal commutator term $[L,Lh_x]h$, the local well-posedness, wave breaking, and ill-posedness of this model in Besov spaces have not been fully studied yet. Inspired by the results mentioned above, we consider the following Cauchy problem:
\begin{align}\label{cauchy}
\left\{\begin{array}{ll}
    h_t+\left(\frac{3}{2}h-\frac{1}{2}\right)h_x
    =\frac{1}{2}Lh_x-\frac{1}{2}[L,Lh_x]h,
    \quad &t\in\mathbb{R}^{+},\ x\in\mathbb{R},\\
    h|_{t=0}=h_0,\quad &x\in\mathbb{R}.
\end{array}\right.
\end{align}
Our main results are stated as follows.
\begin{theorem}[\normalfont\itshape {Local well-posedness}]\label{LW}
Let $s\in \mathbb{R}$, $1\leq p,r\leq\infty$ and $(s,p,r)$ satisfy
\begin{align*}
s>1+\frac{1}{p},\quad 1\leq p,r\leq\infty
\quad\text{or}\quad
s=1+\frac{1}{p},\quad 1\leq p<\infty,\quad r=1.
\end{align*}
Assume that $h_0\in B_{p,r}^s(\mathbb{R})$, and set
\begin{align*}
E_{p,r}^{s}(T)\triangleq
\left\{\begin{array}{l}
C([0,T];B_{p,r}^{s}) \cap C^{1}([0,T];B_{p,r}^{s-1}) \quad \text{if} \quad r < \infty, \\
C_{w}([0,T];B_{p,\infty}^{s}) \cap C^{0,1}([0,T];B_{p,\infty}^{s-1})\quad \text{if} \quad r =\infty.
\end{array}\right.
\end{align*}
Then there exists a time $T>0$ such that the Cauchy problem \eqref{cauchy} has a unique solution $h\in E_{p,r}^s(T)$. Moreover, the map $h_0 \mapsto h$ is continuous from $B_{p,r}^s(\mathbb{R})$ into
$C([0,T]; B_{p,r}^{s'}(\mathbb{R}))\cap C^1([0,T]; B_{p,r}^{s'-1}(\mathbb{R}))$
for all $s'<s$ when $r=\infty$ and $s'=s$ when $1\leq r<\infty$.
\end{theorem}
\begin{remark}
\normalfont
Indeed, when $p=r=2$, the Besov space $B^s_{2,2}(\mathbb{R})$ coincides with $H^s(\mathbb{R})$. Therefore, our result covers the well-posedness result of Alonso-Or\'{a}n et al.~\cite{r1}, in which local well-posedness in $H^s(\mathbb{R})$ for $s>\frac{3}{2}$ was established using a priori energy estimates and an approximation procedure using mollifiers.
\end{remark}

\begin{theorem}[\normalfont\itshape{Wave breaking}]\label{th.4}
Let $h_0\in H^s(\mathbb{R})$ with $s>3/2$, and the maximal existence time $T^*>0$ be the lifespan of Cauchy problem \eqref{cauchy}. Fix some $T_0 \in \left(0,\,2 \ln \left( 1 + \frac{1}{2\|h_0\|_{L^\infty}} \right)\right)$
and suppose that there exists a point $x_0 \in \mathbb{R}$ such that
\begin{align}\label{czjs}
h_{0,x}(x_0)\leq -\frac{2\sqrt{\frac{2K}{3}}}{\mathrm{e}^{\sqrt{6K}\,T_0}-1}-\sqrt{\frac{2K}{3}},
\end{align}
where
\begin{align*}
K=2\left(\frac{1}{2 - \frac{4\|h_0\|_{L^\infty}}{2\|h_0\|_{L^\infty} + 1} \mathrm{e}^{\frac{T_0}{2}}} - \frac{1}{4}  \right)^2-\frac{1}{8}.
\end{align*}
Then the solution $h$ blows up in finite time with $T^* \le T_0$.
More precisely,
\begin{align*}
 T^*\leq\frac{1}{\sqrt{6K}}\ln\left(\frac{h_{0,x}(x_0)-\sqrt{\frac{2K}{3}}}{h_{0,x}(x_0)+\sqrt{\frac{2K}{3}}}
\right)
\leq T_0.
\end{align*}
\end{theorem}
\begin{remark}
\normalfont
In contrast to the previous results
\cite{r1,r31,r32,r2}, the wave breaking criterion established in Theorem
\ref{th.4} is obtained without using the $L^2$-conservation law. Instead, the proof is based only on the structural properties of Cauchy problem \eqref{cauchy}.
\end{remark}

\begin{theorem}[\normalfont\itshape{Ill-posedness}]\label{ill}
Let $1\leq p\leq \infty$ and $1< r\leq \infty$. Then for any $N\in\mathbb{N}^+$ large enough, there exists an $h_0\in H^{\infty}(\mathbb{R})$ such that the following hold:\\
$(1)$$~\|h_0\|_{B^{1+\frac 1 p}_{p,r}} \leq \frac{C}{\ln N}\to 0$, as $N\to\infty$;\\
$(2)~$There exists a unique solution $h\in C([0,T);H^{\infty}(\mathbb{R}))$ to the Cauchy problem \eqref{cauchy} with a maximal lifespan
$T\leq \frac{1}{\ln N}$;\\
$(3)~$
$\|h\|_{L^\infty\left([0,T);B^{1+\frac 1 p}_{p,r}\right)}\geq \ln N$.
\end{theorem}
\begin{remark}
\normalfont
For $1\leq p<\infty$, Theorem~\ref{LW} establishes local well-posedness of the Cauchy problem \eqref{cauchy} in $B_{p,1}^{1+\frac{1}{p}}(\mathbb{R})$. This fact, combined with our norm-inflation result in $B_{p,r}^{1+\frac{1}{p}}(\mathbb{R})$ for $1<r\leq\infty$, shows that the restriction $r=1$ is sharp at critical regularity. In particular, taking $p=r=2$, our norm-inflation result implies ill-posedness in the critical Sobolev space $H^{\frac{3}{2}}(\mathbb{R})$. Thus, the local well-posedness result in $H^s(\mathbb{R})$ for $s>\frac{3}{2}$ is also sharp.
\end{remark}

\noindent\textbf{Organization and Strategy}\quad This paper is organized as follows. In Section \ref{sec2}, we review Littlewood-Paley theory, the definition of Besov spaces, and provide some preliminary results and lemmas. In Section \ref{sec}, we are devoted to proving our main results in this paper.

The proof of local well-posedness is given in
Section~\ref{sec3}. The main difficulty in analyzing Eq.~\eqref{dxfc} lies in handling the commutator term, which introduces a more intricate nonlocal structure than that of the FW equation \eqref{eq8} and requires suitable Besov estimates for the approximate solutions
and their differences. To this end, following the ideas in \cite{r51,r59,r202,r201}, we construct approximate solutions by solving a sequence of linear transport equations. Littlewood-Paley theory and transport
estimates then yield uniform bounds in $B^s_{p,r}$ and convergence in $B^{s-1}_{p,r}$. The difference estimates also yield uniqueness, while continuous dependence for $r<\infty$ is established by decomposing the spatial derivative and applying stability estimates
for linear transport equations. These arguments
establish local well-posedness for the Cauchy problem \eqref{cauchy} in $B_{p,r}^{s}(\mathbb{R})$ for $s>1+\frac{1}{p}$ and $1\leq p,r\leq\infty$, as well as at the critical index $s=1+\frac{1}{p}$ for $1\leq p<\infty$ and $r=1$, thereby extending the Sobolev well-posedness theory in \cite{r1} to Besov spaces.

In Section~\ref{sec4}, we explore wave breaking. Our main strategy involves analyzing the evolution of the slope of the solution along characteristic curves. A key step in this analysis is to control the $L^\infty$-norm of strong solutions, but the commutator term makes it difficult to obtain such control directly.
To overcome this difficulty, we make essential use of Osgood's inequality to derive a priori estimates for this norm. Based on these estimates and without using any conservation
law, we establish a new sufficient condition for wave breaking, which differs from those in \cite{r1,r31,r32,r2}.

Finally, motivated by the wave breaking result,
we prove ill-posedness in Section~\ref{sec5}.
Relying on the sufficient condition for wave breaking and the upper bound on the lifespan provided in Theorem~\ref{th.4}, we construct a sequence of smooth initial data with arbitrarily small critical Besov norms and sufficiently negative slopes, so that the corresponding solutions exhibit wave breaking in arbitrarily short time. Combining this construction with the critical blow-up criterion and the embedding $B_{p,r}^{1+\frac{1}{p}}(\mathbb{R}) \hookrightarrow B_{\infty,\infty}^{1}(\mathbb{R})$, we obtain norm inflation, thereby establishing the ill-posedness of the Cauchy problem \eqref{cauchy}
in $B_{p,r}^{1+\frac{1}{p}}(\mathbb{R})$ for
$1\leq p\leq\infty$ and $1<r\leq\infty$.

\section{Preliminaries}\label{sec2}
In this section, we first present some propositions about the Littlewood-Paley decomposition and the nonhomogeneous Besov spaces with their properties.

\begin{proposition}[Littlewood-Paley decomposition \cite{BCD}]
There exists a pair of smooth radial functions $(\chi,\varphi)$
valued in $[0,1]$, such that $\chi$ is supported in the ball
$\mathcal{B}\triangleq\{\xi\in\mathbb{R}^d:|\xi|\leq\frac43\}$ and $\varphi$ is
supported in the annulus
$\mathcal{C}\triangleq\{\xi\in\mathbb{R}^d:
\frac34\leq|\xi|\leq\frac83\}$.
Moreover,
\[
\chi(\xi)+\sum_{j\geq0}\varphi(2^{-j}\xi)=1,
\qquad \forall\,\xi\in\mathbb{R}^d,
\]
and, for all integers $j,j'\geq0$,
\begin{align*}
\operatorname{supp}\varphi(2^{-j}\cdot)
\cap\operatorname{supp}\varphi(2^{-j'}\cdot)
&=\emptyset,
\quad\text{if }|j-j'|\geq2,\\
\operatorname{supp}\chi
\cap\operatorname{supp}\varphi(2^{-j}\cdot)
&=\emptyset,
\quad\text{if }j\geq1.
\end{align*}
\end{proposition}

Then for all $u\in\mathcal{S}^\prime$, we can define the nonhomogeneous dyadic blocks as follows. Let
$$
\Delta_j u\triangleq0,\ if\ j\leq-2,\ \Delta_{-1}u\triangleq\chi(D)u=\mathcal{F}^{-1}(\chi\mathcal{F}u),$$
$$
\Delta_j u\triangleq \varphi(2^{-j}D)u
=\mathcal{F}^{-1}\bigl(\varphi(2^{-j}\cdot)\mathcal{F}u\bigr),
  \ if\ j\geq 0.
$$
Hence
$$
u=\sum_{j\in\mathbb{Z}}\Delta_{j}u\ in\ \mathcal{S}^\prime(\mathbb{R}^d),$$
where the right-hand side is called the nonhomogeneous Littlewood-Paley decomposition of $u$.
		
According to Young's inequality, we get
$$
\Vert\Delta_j u\Vert_{L^p},\Vert S_j u\Vert_{L^p}\leq C\Vert u\Vert_{L^p},\ \forall 1\leq p\leq\infty.$$
The low frequency cut-off operator $S_j$ is defined by
$$
S_j u\triangleq \sum_{j^\prime=-1}^{j-1}\Delta_{j^\prime}u=\chi(2^{-j}D)u=\mathcal{F}^{-1}(\chi(2^{-j}\xi)\mathcal{F}u),\ \forall j\in\mathbb{N}$$
and the Littlewood-Paley decomposition is quasi-orthogonal in $L^2$ in the following sense:
$$
\Delta_{j}\Delta_{k}u\equiv0,\ if\ |j-k|\geq 2,
$$
$$
\Delta_{j}(S_{k-1}u\Delta_{k}u)\equiv0,\ if\ |j-k|\geq 5.$$
		
\begin{definition}[Nonhomogeneous Besov space \cite{BCD}]
Let $s\in\mathbb{R}$, $1\leq p,\ r\leq\infty$. The nonhomogeneous Besov space $B^s_{p,r}(\mathbb{R}^d)${\rm(}$B^s_{p,r}$ for short{\rm)} is defined by
\begin{align*}
B^s_{p,r}(\mathbb{R}^d)\triangleq\left\{f\in\mathcal{S}^\prime:\Vert f\Vert_{B^s_{p,r}}\triangleq\Vert(2^{js}\Vert\Delta_j f\Vert_{L^p})_{j\geq -1}\Vert_{l^r}<\infty\right\}.
 \end{align*}
If $s=\infty$, we have
\begin{align*} B^\infty_{p,r}\triangleq\underset{{s\in\mathbb{R}}}{\cap}B^s_{p,r}=\left\{f\in\mathcal{S}^\prime:\sup_{j\geq -1}2^{js}\Vert\Delta_j f\Vert_{L^p}<\infty\right\}.
\end{align*}
\end{definition}

In the following lemma, we list some important properties of Besov spaces.
\begin{lemma}[\cite{BCD}]\label{Besov}
Let $s\in\mathbb{R},\ 1\leq p,p_1,p_2,r,r_1,r_2\leq\infty$. We have  \\
{\rm(1)} $B^s_{p,r}$ is a Banach space, and is continuously embedded in $\mathcal{S}'$. \\
{\rm(2)} If $r<\infty$, then $\lim\limits_{j\rightarrow\infty}\|S_j u-u\|_{B^s_{p,r}}=0$. If $p,r<\infty$, $C_0^{\infty}$ is dense in $B^s_{p,r}$. \\
{\rm(3)} If $p_1\leq p_2$ and $r_1\leq r_2$, then $ B^s_{p_1,r_1}\hookrightarrow B^{s-d(\frac 1 {p_1}-\frac 1 {p_2})}_{p_2,r_2}$.
If $s_1<s_2$, the embedding $B^{s_2}_{p,r_2}\hookrightarrow B^{s_1}_{p,r_1}$ is locally compact. \\
{\rm(4)} $\forall s>0$, $B^s_{p,r}\cap L^\infty$ is an algebra. Moreover, $B^s_{p,r}\hookrightarrow L^{\infty} \Leftrightarrow s>\frac d p\ \text{or}\ s=\frac d p,\ r=1$ $\Leftrightarrow B^s_{p,r}$ is an algebra. \\
{\rm(5)} Complex interpolation: $\forall f\in B^{s_1}_{p,r}\cap B^{s_2}_{p,r},\ \forall\theta\in[0,1]$
\begin{align}
	\Vert f\Vert_{B^{\theta s_1+(1-\theta) s_2}_{p,r}}\leq\Vert f\Vert^\theta_{B^{s_1}_{p,r}}\Vert f\Vert^{1-\theta}_{B^{s_2}_{p,r}}.
\end{align}
{\rm(6)} Logarithmic interpolation: $\forall s\in\mathbb{R},\ \varepsilon>0$, and $1\leq p\leq\infty$, there exists a constant $C$ such that
\begin{align}
\Vert u\Vert_{B^s_{p,1}}\leq C\frac{1+\varepsilon} {\varepsilon}\Vert u\Vert_{B^s_{p,\infty}}\left(1+\log\frac{\Vert u\Vert_{B^{s+\varepsilon}_{p,\infty}}}{\Vert u\Vert_{B^{s}_{p,\infty}}}\right).
\end{align}\\
{\rm(7)} Fatou property: if $(u_n)_{n\in\mathbb{N}}$ is a bounded sequence in $B^s_{p,r}$, then an element $u\in B^s_{p,r}$ and a subsequence $(u_{n_k})_{k\in\mathbb{N}}$ exist such that
$$ \lim_{k\rightarrow\infty}u_{n_k}=u\ \text{in}\ \mathcal{S}'\quad \text{and}\quad \|u\|_{B^s_{p,r}}\leq C\liminf_{k\rightarrow\infty}\|u_{n_k}\|_{B^s_{p,r}}. $$
{\rm(8)} Let $m\in\mathbb{R}$ and $f$ be an $S^m$- multiplier, (i.e. $f$ is a smooth function and satisfies that $\forall\alpha\in\mathbb{N}^d$, $\exists C=C(\alpha)$, such that $|\partial^{\alpha}f(\xi)|\leq C(1+|\xi|)^{m-|\alpha|},\ \forall\xi\in\mathbb{R}^d)$.
Then the operator $f(D)=\mathcal{F}^{-1}(f\mathcal{F}\cdot)$ is continuous from $B^s_{p,r}$ to $B^{s-m}_{p,r}$.
\end{lemma}

The 1-D Moser-type estimates are provided as follows.
\begin{lemma}[\cite{BCD}]\label{product}
The following estimates hold:\\
(1) For any $s>0$ and any $(p,r)$ in $[1,\infty]^2$, the space $L^{\infty} \cap B^s_{p,r}$ is an algebra, and a constant $C=C(s)$ exists such that
\begin{align*}
\|uv\|_{B^s_{p,r}}\leq C(\|u\|_{L^{\infty}}\|v\|_{B^s_{p,r}}+\|u\|_{B^s_{p,r}}\|v\|_{L^{\infty}}). 
\end{align*}
(2) If $1\leq p,r\leq \infty,\ s_1\leq s_2,\ s_2>\frac{1}{p}\ (s_2 \geq \frac{1}{p}\ \text{if}\ r=1)$ and $s_1+s_2>\max\{0, \frac{2}{p}-1\}$, there exists $C=C(s_1,s_2,p,r)$ such that
\begin{align*}
\|uv\|_{B^{s_1}_{p,r}}\leq C\|u\|_{B^{s_1}_{p,r}}\|v\|_{B^{s_2}_{p,r}}. \end{align*}
\end{lemma}

Next, we state some useful results in the transport equation theory, which are crucial to the proofs of Theorem \ref{LW}. Let us consider the linear transport equation:
	\begin{align}
		\left\{
		\begin{array}{l}\label{transport}
			\partial_t f+v\cdot \nabla f=g,\\
			f|_{t=0}=f_0.
		\end{array}
		\right.
	\end{align}
	
\begin{lemma}[A priori estimates in Besov spaces \cite{BCD}]\label{priori estimate}
Let $1\leq p,r\leq\infty,\ s>-\min\{\frac{1}{p},\frac{1}{p'}\}$.There exists a constant $C$ such that for all solutions $f\in L^{\infty}([0,T];B^s_{p,r})$ of \eqref{transport} with initial data $f_0\in B^s_{p,r}$, and $g\in L^1([0,T];B^s_{p,r})$, we have, for a.e. $t\in[0,T]$, 
\begin{align*} 
\|f(t)\|_{B^s_{p,r}}\leq \|f_0\|_{B^s_{p,r}}+\int_0^t \|g(t')\|_{B^s_{p,r}}dt'+C\int_0^t V^{\prime} (t')\|f(t')\|_{B^s_{p,r}}dt{'} 
\end{align*}
or
\begin{align*} 
\|f(t)\|_{B^s_{p,r}}\leq e^{CV(t)}\Big(\|f_0\|_{B^s_{p,r}}+\int_0^t e^{-CV(t')}\|g(t')\|_{B^s_{p,r}}dt'\Big) 
\end{align*}
with
\begin{align*}
V^{\prime}(t)=\left\{
\begin{array}{ll}
     \|\partial_xv(t)\|_{B_{p,\infty}^{\frac{1}{p}}\cap L^\infty}, & \text{if} \ s<1+\frac{1}{p}, \\
     \|\partial_xv(t)\|_{B_{p,r}^s}, & \text{if}\ s=1+\frac{1}{p}\ and\ r>1, \\
     \|\partial_xv(t)\|_{B_{p,r}^{s-1}}, & \text{if}\ s>1+\frac{1}{p}\ or\ \left\{s=1+\frac{1}{p}\ and\ r=1\right\}.
\end{array}\right.
\end{align*}
In particular, if $f=v$, then for all $s>0$, $V^{\prime}(t)=\|\partial_x  v(t)\|_{L^{\infty}}.$
\end{lemma}
\begin{lemma}[\cite{BCD}]\label{existence}
Let $1\leq p, r\leq\infty,\ s> -\min\{\frac 1 p, \frac 1 {p'}\}$. Suppose $f_0\in B^s_{p,r}$, $g\in L^1([0,T];B^s_{p,r})$. Let $v$ be a time-dependent vector field such that $v\in L^\rho([0,T];B^{-M}_{\infty,\infty})$ for some $\rho>1$ and $M>0$, and
\begin{align*}
\begin{array}{ll}
	\partial_x v\in L^1([0,T];B^{\frac 1 p}_{p,\infty}\cap L^\infty), &\ \text{if}\ s<1+\frac 1 p, \\
	\partial_x v\in L^1([0,T];B^{s}_{p,r}), &\ \text{if}\ s=1+\frac 1 p,\ r>1, \\
	\partial_x v\in L^1([0,T];B^{s-1}_{p,r}), &\ \text{if}\ s>1+\frac 1 {p}\ or\ \left\{s=1+\frac{1}{p}\ and\ r=1\right\}.
\end{array}
\end{align*}
Then the equation \eqref{transport} has a unique solution $f$ in
\begin{itemize}
\item[-] the space $C([0,T];B^s_{p,r})$, if $r<\infty$;
\item[-] the space
$\Big(\underset{{s'<s}}{\bigcap} C([0,T];B^{s'}_{p,\infty})\Big)
\cap C_w([0,T];B^s_{p,\infty})$, if $r=\infty$.
\end{itemize}
\end{lemma}

\begin{lemma}[\cite{BCD}]\label{continuous}
Let $1\leq p\leq\infty,\ 1\leq r<\infty,\ s>1+\frac 1 p\ (or\ s=1+\frac 1 p,\ 1\leq p<\infty,\ r=1)$. Let $\overline{\mathbb{N}}=\mathbb{N}\cup\{\infty\}$. For $n\in\overline{\mathbb{N}}$, assume that $f^n\in C([0,T];B^{s-1}_{p,r})$ is the solution to
\begin{align*}
   	\left\{\begin{array}{l}
   		\partial_t f^n+a^n\partial_x f^n=g, \\
   		f^n|_{t=0}=f_0
 	\end{array}\right.
\end{align*}
with $f_0\in B^{s-1}_{p,r},\ g\in L^1([0,T];B^{s-1}_{p,r})$. In addition, suppose that for some $\alpha\in L^1([0,T])$,
\begin{align*}
\sup\limits_{n\in\overline{\mathbb{N}}}\|a^n(t)\|_{B^{s}_{p,r}}\leq \alpha(t).
\end{align*}
If $a^n$ converges to $a^{\infty}$ in $L^1([0,T];B^{s-1}_{p,r})$, then $f^n$ converges to $f^{\infty}$ in $C([0,T];B^{s-1}_{p,r})$.
\end{lemma}

We also need the following Kato-Ponce commutator estimate.
\begin{lemma}[\cite{r404}]\label{KP}
Let $s>0$, then there exists a constant $C>0$ such that
\begin{align*}
\|[\Lambda^s,f]g\|_{L^2}
\leq
C\Big(
\|\partial_x f\|_{L^\infty}\|\Lambda^{s-1}g\|_{L^2}
+\|\Lambda^s f\|_{L^2}\|g\|_{L^\infty}
\Big),
\end{align*}
where $\Lambda^s=(1-\partial_x^2)^{s/2}$.
\end{lemma}

Then, the blow-up criteria for Eq.~\eqref{dxfc} can be stated as follows.
\begin{lemma}[\cite{r1}]\label{L2}
Let $h_0\in H^s(\mathbb{R})$ with $s>3/2$, and let maximal existence time $T>0$ be the lifespan associated to the solution $h$ to  \eqref{dxfc} with $h(0,x)=h_0(x)$, then $h(t,x)$ blows up in finite time $T$ if and only if
\begin{align*}
\int_0^{T}\|h_x(\tau)\|_{BMO}d\tau=\infty.
\end{align*}
\end{lemma}

In addition, the following lemma shows that for the maximal time of existence $T> 0$, the solution remains bounded.
\begin{lemma}[\cite{r1}]\label{L3}
Let $h_0\in H^s(\mathbb{R})$ with $s>3/2$, and let $T>0$ be the lifespan associated to the solution $h$ to  \eqref{dxfc} with $h(0,x)=h_0(x)$, then
\begin{align*}
\sup_{t\in [0,T) }\|h(t)\|_{L^\infty(\mathbb{R})}<\infty.
\end{align*}
\end{lemma}

Finally, we recall the following key lemma, which will be used in the proof of our wave breaking results.

\begin{lemma}[Osgood Inequality \cite{BCD}]\label{Osgood}
Let $f$ be a positive measurable function, $\lambda$ a positive locally integrable function and $\mu$ a positive increasing continuous function. If
\begin{align*}
    f(t) \leq a + \int_{t_0}^{t} \lambda(s)\mu(f(s)) \, ds, \quad \text{for } a > 0,
\end{align*}
then
\begin{align*}
    W(f(t)) \leq W(a) + \int_{t_0}^{t} \lambda(s) \, ds,
\end{align*}
where $W(x) := \int_{b}^{x} \frac{d\tau}{\mu(\tau)}$ for some $b > 0$.
\end{lemma}

\section{Proofs of main results}\label{sec}
\subsection{The proof of local well-posedness}\label{sec3}
In this section, we establish the local well-posedness in Besov spaces for the Cauchy problem \eqref{cauchy}. 
\begin{proof}[\normalfont\bfseries Proof of Theorem~\ref{LW}]
Now  we prove Theorem \ref{LW} in the following steps.\\

\noindent\textbf{Step 1. Constructing Approximate Solutions and Uniform Bounds.}

Let $h^0\triangleq0$ and define a sequence $(h^n)_{n\in\mathbb{N}}$ of smooth functions by solving the following linear transport equation
\begin{align}\label{chushi}
\left\{\begin{array}{l}
\partial_t h^{n+1}
    +\left(\frac{3}{2}h^n-\frac{1}{2}\right)\partial_xh^{n+1}
    =\frac{1}{2}L\partial_xh^n
    -\frac{1}{2}[L,L\partial_xh^n]h^n,\\
h^{n+1}|_{t=0}=S_{n+1}h_{0}.
\end{array}\right.
\end{align}
Since $h^0=0$, we have $h^1(t,x)=(S_1h_0)(x+\frac{t}{2})$ and hence $h^1\in E_{p,r}^{s}(T)$ for all $T>0$. By induction, we assume that $h^n\in E_{p,r}^{s}(T)$ for all $T>0$. Since $L=(1-\partial_x^2)^{-1}$ is an $S^{-2}$-multiplier and $L\partial_x$ is an $S^{-1}$-multiplier, we have 
\begin{align*}
\|L\partial_x h^n\|_{B^s_{p,r}}
\leq C\|h^n\|_{B^{s-1}_{p,r}}
\leq C\|h^n\|_{B^s_{p,r}}.
\end{align*}
Moreover, under the assumptions of Theorem \ref{LW}, the algebra properties of $B^s_{p,r}$ and $B^{s-1}_{p,r}$, together with the embedding $B^{s-1}_{p,r}\hookrightarrow L^\infty$, yield
\begin{align*}
\big\|[L,L\partial_xh^n]h^n\big\|_{B^s_{p,r}}
&=\big\|L\big((L\partial_xh^n)h^n\big)-
(L\partial_xh^n)Lh^n\big\|_{B^s_{p,r}}\\
&\leq C\|L\partial_xh^n\|_{B^s_{p,r}}\|h^n\|_{B^s_{p,r}}+C\|L\partial_xh^n\|_{B^s_{p,r}}\|Lh^n\|_{B^s_{p,r}}\\
&\leq C\|h^n\|_{B^s_{p,r}}^2.
\end{align*}
Therefore,
\begin{align}\label{yxgj}
\left\|
\frac{1}{2}L\partial_xh^n
-\frac{1}{2}[L,L\partial_xh^n]h^n
\right\|_{B^s_{p,r}}
\leq C\left(\|h^n\|_{B^s_{p,r}}
+\|h^n\|_{B^s_{p,r}}^2\right).
\end{align}
Then by Lemma\ \ref{Besov} and Lemma\ \ref{product} , we have
\begin{align*}
\frac{1}{2}L\partial_xh^n-\frac{1}{2}[L,L\partial_xh^n]h^n
\in L^1([0,T];B^s_{p,r}),\ \forall\ T>0.
\end{align*}
Applying Lemma \ref{existence} to \eqref{chushi}, we deduce that \eqref{chushi} has a solution $h^{n+1}\in E_{p,r}^{s}(T),\ \forall\ T>0$.

Define $H^n(t)=\int_0^t\|h^{n}(t')\|_{B^{s}_{p,r}}dt'$. Since $\|h^{n+1}(0)\|_{B^{s}_{p,r}} =\|S_{n+1}h_{0}\|_{B^{s}_{p,r}} \leq C\|h_0\|_{B^{s}_{p,r}} $, by Lemma \ref{priori estimate} and \eqref{yxgj}, we have
\begin{align}\label{LC1}
\|h^{n+1}(t)\|_{B^{s}_{p,r}}
&\leq C\|h_0\|_{B^{s}_{p,r}}\mathrm{e}^{CH^n(t)}+\int_0^t \mathrm{e}^{C\int_{t'}^t \|h^{n}(t'')\|_{B^{s}_{p,r}}dt''} \left\|\left(\frac{1}{2}L\partial_xh^n-\frac{1}{2}[L,L\partial_xh^n]h^n\right)(t')\right\|_{B^{s}_{p,r}}dt'\nonumber\\
&\leq C\Big(\|h_0\|_{B^{s}_{p,r}}\mathrm{e}^{CH^n(t)}+\int_0^t \mathrm{e}^{C\int_{t'}^t \|h^{n}(t'')\|_{B^{s}_{p,r}}dt''}\left(\|h^{n}(t')\|_{B^{s}_{p,r}}+\|h^{n}(t')\|_{B^{s}_{p,r}}^2\right)dt'\Big).
\end{align}
For fixed $T>0$ such that
$2C^2T\big(1+\|h_0\|_{B^{s}_{p,r}}\big)<1$, we have
\begin{align*}
\|h^1(t)\|_{B^{s}_{p,r}}
=\|S_1h_0\|_{B^{s}_{p,r}}
\leq\frac{C\|h_0\|_{B^{s}_{p,r}}}
{1-2C^2t\big(1+\|h_0\|_{B^{s}_{p,r}}\big)},
\ \forall\ t\in[0,T],
\end{align*}
By induction, we assume that
\begin{align}\label{LC2}
\|h^n(t)\|_{B^{s}_{p,r}}
\leq\frac{C\|h_0\|_{B^{s}_{p,r}}}
{1-2C^2t\big(1+\|h_0\|_{B^{s}_{p,r}}\big)},
\ \forall\ t\in[0,T].
\end{align}
Substituting \eqref{LC2} into \eqref{LC1}, we have
\begin{align*}
\|h^{n+1}(t)\|_{B^{s}_{p,r}}
&\leq C\|h_0\|_{B^{s}_{p,r}}
\Big(1-2C^2t\big(1+\|h_0\|_{B^{s}_{p,r}}\big)
\Big)^{-\frac12}+C\int_0^t
\Bigg(\frac{1-2C^2t\big(1+\|h_0\|_{B^{s}_{p,r}}\big)}{
1-2C^2t'\big(1+\|h_0\|_{B^{s}_{p,r}}\big)}
\Bigg)^{-\frac12}\\
&\qquad\qquad\times
\left[\frac{C\|h_0\|_{B^{s}_{p,r}}}
{1-2C^2t'\big(1+\|h_0\|_{B^{s}_{p,r}}\big)}
+\Bigg(\frac{C\|h_0\|_{B^{s}_{p,r}}}
{1-2C^2t'\big(1+\|h_0\|_{B^{s}_{p,r}}\big)}
\Bigg)^2\right]dt'\\
&\leq C\|h_0\|_{B^{s}_{p,r}}
\Big(1-2C^2t\big(1+\|h_0\|_{B^{s}_{p,r}}\big)
\Big)^{-\frac12}+C^3\|h_0\|_{B^{s}_{p,r}}
\big(1+\|h_0\|_{B^{s}_{p,r}}\big)\\
&\qquad\qquad\times\int_0^t\Bigg(\frac{1-2C^2t\big(1+\|h_0\|_{B^{s}_{p,r}}\big)
}{1-2C^2t'\big(1+\|h_0\|_{B^{s}_{p,r}}\big)
}\Bigg)^{-\frac12}
\Big(1-2C^2t'\big(1+\|h_0\|_{B^{s}_{p,r}}\big)
\Big)^{-2}dt'\\
&\leq\frac{C\|h_0\|_{B^{s}_{p,r}}}
{1-2C^2t\big(1+\|h_0\|_{B^{s}_{p,r}}\big)},
\ \forall\ t\in[0,T].
\end{align*}
Therefore, $(h^n)_{n\in \mathbb{N}}$ is uniformly bounded in $L^{\infty}([0,T];B^{s}_{p,r})$, which implies that $\left(
\left(\frac{3}{2}h^n-\frac{1}{2}\right)
\partial_x h^{n+1}
\right)_{n\in\mathbb{N}}$ is uniformly  bounded in $L^{\infty}([0,T];B^{s-1}_{p,r})$ and the right hand side of \eqref{chushi} is uniformly  bounded in $L^{\infty}([0,T];B^{s}_{p,r})$. We can conclude that the sequence $(h^n)_{n\in \mathbb{N}}$ is uniformly bounded in $E^{s}_{p,r}(T)$.\\

\noindent\textbf{Step 2. Convergence and Regularity.}

We next show that $(h^n)_{n\in\mathbb{N}}$ is a Cauchy
sequence in $C([0,T];B^{s-1}_{p,r})$. To this end, for any $(m,n)\in\mathbb{N}^2$, we have
\begin{align}\label{cauchy1}
\left\{
\begin{array}{l}
\left(\partial_t+\left(\frac{3}{2}h^{n+m}-\frac{1}{2}\right)\partial_x\right)
(h^{n+m+1}-h^{n+1})
=\frac{3}{2}(h^n-h^{n+m})\partial_xh^{n+1}\\
\qquad\qquad
+\frac{1}{2}\Big(
L\partial_x(h^{n+m}-h^n)
-[L,L\partial_xh^{n+m}]h^{n+m}
+[L,L\partial_xh^n]h^n
\Big)
\triangleq g^{n,m},\\
(h^{n+m+1}-h^{n+1})|_{t=0}
=S_{n+m+1}h_{0}-S_{n+1}h_{0}.
\end{array}\right.
\end{align}
Define
\begin{align*}
\widetilde{H}^{n,m}(t)=
\left\{
\begin{array}{ll}
\|\partial_x h^{n+m}(t)\|_{B_{p,\infty}^{\frac{1}{p}}\cap L^\infty},
& \text{if } s<2+\frac{1}{p}, \\[1mm]
\|\partial_x h^{n+m}(t)\|_{B_{p,r}^{s-1}},
& \text{if } s=2+\frac{1}{p},\ r>1, \\[1mm]
\|\partial_x h^{n+m}(t)\|_{B_{p,r}^{s-2}},
& \text{if } s>2+\frac{1}{p}
\ \text{or } \left(s=2+\frac{1}{p},\ r=1\right).
\end{array}
\right.
\end{align*}
Since $(h^n)_{n\in \mathbb{N}}$ is uniformly bounded in $L^{\infty}([0,T];B^{s}_{p,r})$, we have $\widetilde{H}^{n,m}(t)\leq C_T,\ \forall\ (m,n)\in\mathbb{N}^2,\ \forall\ t\in [0,T]$. Applying Lemma \ref{priori estimate} to \eqref{cauchy1}, we have
\begin{align*}
&\|(h^{n+m+1}-h^{n+1})(t)\|_{B^{s-1}_{p,r}}\\
\leq&\ e^{C\int_0^t \widetilde{H}^{n,m}(t')dt'}
\Bigg(\|S_{n+m+1}h_{0}-S_{n+1}h_{0}\|_{B^{s-1}_{p,r}}
+\int_0^t e^{-C\int_0^{t'}\widetilde{H}^{n,m}(t'')dt''}\|g^{n,m}(t')\|_{B^{s-1}_{p,r}}dt'
\Bigg)\\
\leq&\ C_T\Bigg(\|S_{n+m+1}h_{0}-S_{n+1}h_{0}\|_{B^{s-1}_{p,r}}+\int_0^t
\|g^{n,m}(t')\|_{B^{s-1}_{p,r}}dt'\Bigg).
\end{align*}
Taking advantage of Lemma \ref{Besov} and Lemma \ref{product}, we  obtain
\begin{align*}
&\|(h^n-h^{n+m})\partial_xh^{n+1}\|_{B^{s-1}_{p,r}}\leq C\|h^n-h^{n+m}\|_{B^{s-1}_{p,r}}
\|\partial_xh^{n+1}\|_{B^{s-1}_{p,r}}\leq C_T\|h^n-h^{n+m}\|_{B^{s-1}_{p,r}},\\
&\|L\partial_x(h^{n+m}-h^n)\|_{B^{s-1}_{p,r}}
\leq C\|h^{n+m}-h^n\|_{B^{s-2}_{p,r}}
\leq C_T\|h^{n+m}-h^n\|_{B^{s-1}_{p,r}},\\
&\|[L,L\partial_xh^n]h^n-[L,L\partial_xh^{n+m}]h^{n+m}\|_{B^{s-1}_{p,r}}
\leq C\|[L,L\partial_x(h^n-h^{n+m})]h^n\|_{B^{s-1}_{p,r}}+C\|[L,L\partial_xh^{n+m}](h^n-h^{n+m})\|_{B^{s-1}_{p,r}}\\
&\qquad\qquad\qquad\qquad\qquad\qquad\qquad\qquad\qquad\leq C\big(\|h^n\|_{B^s_{p,r}}+\|h^{n+m}\|_{B^s_{p,r}}\big)\|h^n-h^{n+m}\|_{B^{s-1}_{p,r}}\\
&\qquad\qquad\qquad\qquad\qquad\qquad\qquad\qquad\qquad\leq C_T\|h^n-h^{n+m}\|_{B^{s-1}_{p,r}}.
\end{align*}
Since $\|S_{n+m+1}h_{0}-S_{n+1}h_{0}\|_{B^{s-1}_{p,r}}\leq C2^{-n}\|h_0\|_{B^s_{p,r}}$, we have
\begin{align*}
&\|(h^{n+m+1}-h^{n+1})(t)\|_{B^{s-1}_{p,r}}\\
\leq\ &C_T\Big(2^{-n}\|h_0\|_{B^s_{p,r}}
+\int_0^t\|h^{n+m}(t_1)-h^n(t_1)\|_{B^{s-1}_{p,r}}dt_1\Big)\\
\leq\ &C_T\Big(2^{-n}\|h_0\|_{B^s_{p,r}}
+\int_0^tC_T\Big(2^{-n+1}\|h_0\|_{B^s_{p,r}}
+\int_0^{t_1}\|h^{n+m-1}(t_2)-h^{n-1}(t_2)\|_{B^{s-1}_{p,r}}dt_2\Big)dt_1\Big)\\
\leq\ &\dots\\
\leq\ &\sum_{l=1}^{n+1}C_T^l2^{-n+l-1}\|h_0\|_{B^s_{p,r}}
\frac{T^{l-1}}{(l-1)!}
+\int_0^t\int_0^{t_1}\dots\int_0^{t_n}
C_T^{n+1}\|h^m(t_{n+1})\|_{B^{s-1}_{p,r}}
dt_{n+1}\dots dt_2dt_1\\
\leq\ &C_T2^{-n}\|h_0\|_{B^s_{p,r}}e^{2TC_T}
+\frac{C_T^{n+2}T^{n+1}}{(n+1)!}
\rightarrow0\quad(n\rightarrow\infty).
\end{align*}
Therefore, $(h^n)_{n\in \mathbb{N}}$ is a Cauchy sequence in $C([0,T];B^{s-1}_{p,r})$. Hence, there exists $h\in C([0,T];B^{s-1}_{p,r})$ such that $h^n\to h$ in $C([0,T];B^{s-1}_{p,r})$.
Finally, we have to check that $h\in E_{p,r}^{s}(T)$ and satisfies \eqref{cauchy}. Combining the uniform boundedness of $(h^n)_{n\in \mathbb{N}}$ in $L^{\infty}([0,T];B^{s}_{p,r})$ and the Fatou property, we can obtain that $h\in L^{\infty}([0,T];B^{s}_{p,r})$. Since $(h^n)_{n\in \mathbb{N}}$ converges to $h$ in $C([0,T];B^{s-1}_{p,r})$, an interpolation argument ensures that the convergence holds true in $C([0,T];B^{s'}_{p,r}),\ \forall\ s'<s$. Thus, we can take the limit in \eqref{chushi} and conclude that $h$ is a solution of \eqref{cauchy}. Since $h\in L^{\infty}([0,T];B^{s}_{p,r})$, by Lemma \ref{Besov} and Lemma \ref{product}, we have
\begin{align*}
\frac{1}{2}L\partial_xh-\frac{1}{2}[L,L\partial_xh]h
\in L^{\infty}([0,T];B^{s}_{p,r}).
\end{align*}
Hence, by Lemma \ref{existence}, we see that $h$ is in $C([0,T];B^s_{p,r})$ (resp.,
$\big(\underset{s'<s}{\bigcap}C([0,T];B^{s'}_{p,\infty})\big)
\bigcap C_w([0,T];B^s_{p,\infty})$) if $r<\infty$ (resp., $r=\infty$). Thus, we deduce from \eqref{cauchy} that $h_t$
belongs to $C([0,T];B^{s-1}_{p,r})$ (resp., $L^{\infty}([0,T];B^{s-1}_{p,\infty})$) if $r<\infty$ (resp., $r=\infty$), and conclude that $h\in E_{p,r}^{s}(T)$.\\

\noindent\textbf{Step 3. Uniqueness.}

Suppose $h,\bar h$ are two solutions of \eqref{cauchy} with initial data $h_0$ and $\bar h_0$, respectively. Then
\begin{equation}\label{unique}
\left\{
\begin{array}{l}
\left(\partial_t+\left(\frac{3}{2}h-\frac{1}{2}\right)\partial_x\right)(h-\bar h)
=\frac{3}{2}(\bar h-h)\partial_x\bar h+\frac{1}{2}L\partial_x(h-\bar h)
-\frac{1}{2}[L,L\partial_xh]h+\frac{1}{2}[L,L\partial_x\bar h]\bar h,\\
(h-\bar h)|_{t=0}=h_0-\bar h_0.
\end{array}
\right.
\end{equation}
Define
\begin{align*}
\widetilde{H}(t)=\left\{
\begin{array}{ll}
\|\partial_x h(t)\|_{B_{p,\infty}^{\frac{1}{p}}\cap L^\infty},
& \text{if } s<2+\frac{1}{p}, \\
\|\partial_x h(t)\|_{B_{p,r}^{s-1}},
& \text{if } s=2+\frac{1}{p},\ r>1, \\
\|\partial_x h(t)\|_{B_{p,r}^{s-2}},
& \text{if } s>2+\frac{1}{p}\ \text{or}\ \left(s=2+\frac{1}{p},\ r=1\right).
\end{array}\right.
\end{align*}
Similar to Step 2, applying Lemma \ref{priori estimate} to \eqref{unique} gives
\begin{align*}
&\mathrm{e}^{-C\int_0^t \widetilde{H}(t')dt'}\|(h-\bar h)(t)\|_{B^{s-1}_{p,r}}\\
\leq\ & \|h_{0}-\bar h_{0}\|_{B^{s-1}_{p,r}}
+\int_0^t \mathrm{e}^{-C\int_0^{t'}\widetilde{H}(t'')dt''}
\Big(\|(\bar h-h)\partial_x\bar h\|_{B^{s-1}_{p,r}}
+\|L\partial_x(h-\bar h)\|_{B^{s-1}_{p,r}}
+\|[L,L\partial_xh]h-[L,L\partial_x\bar h]\bar h\|_{B^{s-1}_{p,r}}\Big)dt'\\
\leq\ & \|h_{0}-\bar h_{0}\|_{B^{s-1}_{p,r}}
+C\int_0^t \mathrm{e}^{-C\int_0^{t'}\widetilde{H}(t'')dt''}
\|h-\bar h\|_{B^{s-1}_{p,r}}
\Big(1+\|h\|_{B^{s}_{p,r}}+\|\bar h\|_{B^{s}_{p,r}}\Big)dt'.
\end{align*}
Taking advantage of Gronwall's inequality, we have
\begin{align*}
\mathrm{e}^{-C\int_0^t \widetilde{H}(t')dt'}\|(h-\bar h)(t)\|_{B^{s-1}_{p,r}}\leq\|h_{0}-\bar h_{0}\|_{B^{s-1}_{p,r}}
e^{C\int_0^t\left(1+\|h(t')\|_{B^{s}_{p,r}}+\|\bar h(t')\|_{B^{s}_{p,r}}\right)dt'}.
\end{align*}
Since $\widetilde{H}(t)\leq C\|h(t)\|_{B^{s}_{p,r}}$, we obtain
\begin{align*}
\|(h-\bar h)(t)\|_{B^{s-1}_{p,r}}
\leq
\|h_{0}-\bar h_{0}\|_{B^{s-1}_{p,r}}
e^{C\int_0^t\left(1+\|h(t')\|_{B^{s}_{p,r}}+\|\bar h(t')\|_{B^{s}_{p,r}}\right)dt'}.
\end{align*}
Hence, if $h_0=\bar h_0$, then $h\equiv\bar h$ and we complete the proof of uniqueness.\\

\noindent\textbf{Step 4. The continuous dependence.}

Denote $\overline{\mathbb{N}}=\mathbb{N}\cup\{\infty\}$. Suppose that $h^n$ is the solution
of \eqref{cauchy} with initial data $h_0^n\in B^{s}_{p,r}$ and that
$h_0^n\to h_0^\infty$ in $B^{s}_{p,r}$ as $n\to\infty$.
Since $(h_0^n)_{n\in\overline{\mathbb N}}$ is uniformly bounded in
$B^{s}_{p,r}$, there exists a $T>0$ such that
$h^n\in E_{p,r}^{s}(T)$ for all $n\in\overline{\mathbb N}$.

\textbf{Case 1: $r<\infty$}

Fix $T>0$ such that
\begin{align*}
2C^2T\Big(1+\underset{n\in\overline{\mathbb{N}}}{\sup}
\|h_0^n\|_{B^{s}_{p,r}}\Big)<1.
\end{align*}
Then by Step 1, we have
\begin{align*}
\underset{n\in\overline{\mathbb{N}}}{\sup}
\|h^n\|_{L^\infty_T(B^{s}_{p,r})}
\leq \frac{C\underset{n\in\overline{\mathbb{N}}}{\sup}
\|h_0^n\|_{B^{s}_{p,r}}}{1-2C^2T\Big(1+\underset{n\in\overline{\mathbb{N}}}{\sup}\|h_0^n\|_{B^{s}_{p,r}}\Big)}
\triangleq M.
\end{align*}
Similar to Step 3, we see that $h^n\to h^\infty$ in
$C([0,T];B^{s-1}_{p,r})$ as $n\to\infty$. Hence, we only
need to prove $\partial_xh^n\to\partial_xh^\infty$ in
$C([0,T];B^{s-1}_{p,r})$. For simplicity, let
\begin{align*}
v^n\triangleq\partial_xh^n,\qquad
f^n\triangleq
-\frac{3}{2}(\partial_xh^n)^2
+\frac{1}{2}(L-I)h^n
-\frac{1}{2}\partial_x\big([L,L\partial_xh^n]h^n\big),
\end{align*}
and decompose $v^n$ into $v^n=z^n+w^n$ with $(z^n,w^n)$ satisfying
\begin{align}\label{zn}
\left\{
\begin{aligned}
&\partial_tz^n+\left(\frac{3}{2}h^n-\frac{1}{2}\right)\partial_xz^n
=f^n-f^\infty,\\
&z^n|_{t=0}=\partial_xh_0^n-\partial_xh_0^\infty
\end{aligned}
\right.
\end{align}
and
\begin{align}\label{wn}
\left\{
\begin{aligned}
&\partial_tw^n+\left(\frac{3}{2}h^n-\frac{1}{2}\right)\partial_xw^n
=f^\infty,\\
&w^n|_{t=0}=\partial_xh_0^\infty.
\end{aligned}
\right.
\end{align}
Similar to Step 1 and Step 2, applying Lemma \ref{priori estimate} to \eqref{zn} gives
\begin{align*}
\|z^n(t)\|_{B^{s-1}_{p,r}}
\leq \mathrm{e}^{CMT}\Big(\|\partial_xh_{0}^n-\partial_xh_{0}^\infty\|_{B^{s-1}_{p,r}}+\int_0^t \|f^n(t')-f^{\infty}(t')\|_{B^{s-1}_{p,r}}dt'\Big).
\end{align*}
In addition, Lemma \ref{Besov} and Lemma \ref{product} imply that
\begin{align*}
\|(h^n_x)^2-(h^{\infty}_x)^2\|_{B^{s-1}_{p,r}}
&\leq C\|h^n_x-h^{\infty}_x\|_{B^{s-1}_{p,r}}\|h^n_x+h^{\infty}_x\|_{B^{s-1}_{p,r}}\\
&\leq CM\|h^n_x-h^{\infty}_x\|_{B^{s-1}_{p,r}},\\
\|(L-I)(h^n-h^{\infty})\|_{B^{s-1}_{p,r}}
&\leq C\|h^n-h^{\infty}\|_{B^{s-1}_{p,r}},\\
\|\partial_x\big([L,L\partial_xh^n]h^n-[L,L\partial_xh^{\infty}]h^{\infty}\big)\|_{B^{s-1}_{p,r}}
&\leq C\|[L,L\partial_xh^n]h^n-[L,L\partial_xh^{\infty}]h^{\infty}\|_{B^{s}_{p,r}}\\
&\leq C\big(\|h^n\|_{B^{s}_{p,r}}+\|h^{\infty}\|_{B^{s}_{p,r}}\big)
\|h^n-h^{\infty}\|_{B^{s-1}_{p,r}}\\
&\leq CM\|h^n-h^{\infty}\|_{B^{s-1}_{p,r}}.
\end{align*}
Therefore,
\begin{align*}
\|z^n(t)\|_{B^{s-1}_{p,r}}
\leq C\mathrm{e}^{CMT}\Bigg(\|h_{0}^n-h_{0}^\infty\|_{B^{s}_{p,r}}
+(1+M)\int_0^t\Big(\|h^n(t')-h^{\infty}(t')\|_{B^{s-1}_{p,r}}
+\|h_x^n(t')-h_x^{\infty}(t')\|_{B^{s-1}_{p,r}}\Big)dt'\Bigg).
\end{align*}
We next consider \eqref{wn}. Since $B^{s-1}_{p,r}$ is an algebra and $B^{s-1}_{p,r}\hookrightarrow L^\infty$, we have
\begin{align*}
\|f^{\infty}\|_{B^{s-1}_{p,r}}
&\leq C\|(\partial_x h^{\infty})^2\|_{B^{s-1}_{p,r}}
+C\|(L-I)h^{\infty}\|_{B^{s-1}_{p,r}}
+C\|\partial_x([L,L\partial_xh^{\infty}]h^{\infty})\|_{B^{s-1}_{p,r}}\\
&\leq C\Big(
\|\partial_x h^{\infty}\|^2_{B^{s-1}_{p,r}}
+\|h^{\infty}\|_{B^{s-1}_{p,r}}
+\|h^{\infty}\|^2_{B^{s}_{p,r}}
\Big)\\
&\leq C(M+M^2).
\end{align*}
Recall that $\partial_x h_0^{\infty}\in B^{s-1}_{p,r}$, $\underset{n\in\overline{\mathbb{N}}}{\sup} \|h^n\|_{L^\infty_T(B^{s}_{p,r})}\leq M$, and $h^n\to h^\infty$ in $C([0,T];B^{s-1}_{p,r})$. By Lemma \ref{existence}, we see that $w^n\in C([0,T];B^{s-1}_{p,r})$.
Then we deduce from Lemma \ref{continuous} that
$w^n\to w^{\infty}$ in $C([0,T];B^{s-1}_{p,r})$.
Hence, for all $\varepsilon>0$, there exists a constant $N\in\mathbb{N}$, such that $\forall\ n>N,\ t\in[0,T]$,
\begin{align*}
&\|w^n(t)-w^{\infty}(t)\|_{B^{s-1}_{p,r}}
<\frac{\varepsilon}{4}e^{-C(1+M)Te^{CMT}},\\
&\|h^n(t)-h^{\infty}(t)\|_{B^{s-1}_{p,r}}
<\frac{\varepsilon}{4C(1+M)Te^{CMT}}e^{-C(1+M)Te^{CMT}}.
\end{align*}
Then we have
\begin{align*}
&\|\partial_x h^n(t)-\partial_x h^{\infty}(t)\|_{B^{s-1}_{p,r}}\\
\leq\ &\|z^n(t)\|_{B^{s-1}_{p,r}}+\|w^n(t)-w^{\infty}(t)\|_{B^{s-1}_{p,r}}\\
\leq\ &Ce^{CMT}\Big(\|h_{0}^n-h_{0}^\infty\|_{B^{s}_{p,r}}
+(1+M)\int_0^t\|\partial_x h^n(t')-\partial_x h^{\infty}(t')\|_{B^{s-1}_{p,r}}dt'\Big)
+\frac{\varepsilon}{2}e^{-C(1+M)Te^{CMT}}.
\end{align*}
By Gronwall's inequality, we obtain
$$
\|\partial_x h^n(t)-\partial_x h^{\infty}(t)\|_{B^{s-1}_{p,r}}
\leq\frac{\varepsilon}{2}
+Ce^{CMT}e^{C(1+M)Te^{CMT}}\|h_{0}^n-h_{0}^\infty\|_{B^{s}_{p,r}}.
$$
Since $h_0^n\to h_0^\infty$ in $B^{s}_{p,r}$, we have
$\partial_x h^n\to\partial_x h^\infty$ in
$C([0,T];B^{s-1}_{p,r})$. Therefore,
$$
\|h^n-h^{\infty}\|_{C([0,T];B^{s}_{p,r})}
\leq C\Big(
\|\partial_x h^n-\partial_x h^{\infty}\|_{C([0,T];B^{s-1}_{p,r})}
+\|h^n-h^{\infty}\|_{C([0,T];B^{s-1}_{p,r})}
\Big)\to 0
$$
when $n\to\infty$, which also implies that
\begin{align*}
&\left\|
\frac{1}{2}L\partial_xh^n-\frac{1}{2}[L,L\partial_xh^n]h^n
-\frac{1}{2}L\partial_xh^\infty+\frac{1}{2}[L,L\partial_xh^\infty]h^\infty
\right\|_{B^{s-1}_{p,r}}\\
&\qquad\leq
C\Big(1+\|h^n\|_{B^s_{p,r}}+\|h^\infty\|_{B^s_{p,r}}\Big)
\|h^n-h^\infty\|_{B^s_{p,r}}
\to 0
\end{align*}
when $n\to\infty$. Since
\begin{align*}
\|h^n h^n_x-h^{\infty}h^{\infty}_x\|_{B^{s-1}_{p,r}}
&\leq C\Big(
\|h^n\|_{B^{s}_{p,r}}\|h^n-h^{\infty}\|_{B^{s}_{p,r}}
+\|h^{\infty}\|_{B^{s}_{p,r}}\|h^n-h^{\infty}\|_{B^{s}_{p,r}}
\Big)\to 0
\end{align*}
when $n\to\infty$, we have
$h^n h^n_x\to h^{\infty}h^{\infty}_x$ in
$C([0,T];B^{s-1}_{p,r})$. Moreover,
$\partial_xh^n\to\partial_xh^\infty$ in $C([0,T];B^{s-1}_{p,r})$. Then by \eqref{cauchy}, we see that
$\partial_t h^n\to\partial_t h^\infty$ in
$C([0,T];B^{s-1}_{p,r})$. Hence, we conclude that $h^n\to h^{\infty}$ in $E^s_{p,r}(T)$.\\

\textbf{Case 2: $r=\infty$}

Similar to Case 1, we have
$\underset{n\in\overline{\mathbb{N}}}{\sup}
\|h^n\|_{L^\infty_T(B^{s}_{p,\infty})}\leq M$
and $\|h^n-h^\infty\|_{L^{\infty}([0,T];B^{s-1}_{p,\infty})}\to 0$
as $n\to\infty$. For any fixed $\phi\in B^{-s}_{p',1}$, we write
\begin{align*}
\langle h^n(t)-h^{\infty}(t),\phi\rangle
&=\langle S_j(h^n(t)-h^{\infty}(t)),\phi\rangle
+\langle (Id-S_j)(h^n(t)-h^{\infty}(t)),\phi\rangle\\
&=\langle h^n(t)-h^{\infty}(t),S_j\phi\rangle
+\langle h^n(t)-h^{\infty}(t),(Id-S_j)\phi\rangle.
\end{align*}
Note that
\begin{align*}
\underset{t\in[0,T]}{\sup}
|\langle h^n(t)-h^{\infty}(t),S_j\phi\rangle|
\leq C\|h^n-h^\infty\|_{L^{\infty}([0,T];B^{s-1}_{p,\infty})}
\|S_j\phi\|_{B^{1-s}_{p',1}}
\end{align*}
and
\begin{align*}
\underset{t\in[0,T]}{\sup}
|\langle h^n(t)-h^{\infty}(t),(Id-S_j)\phi\rangle|
&\leq C\|h^n-h^\infty\|_{L^{\infty}([0,T];B^{s}_{p,\infty})}
\|(Id-S_j)\phi\|_{B^{-s}_{p',1}}\\
&\leq 2CM\|(Id-S_j)\phi\|_{B^{-s}_{p',1}}
\to 0\quad (j\to\infty).
\end{align*}
For all $\varepsilon>0$, we can fix $j\in\mathbb{N}$ large enough such that
\begin{align*}
\underset{t\in[0,T]}{\sup}
|\langle h^n(t)-h^{\infty}(t),(Id-S_j)\phi\rangle|
<\frac{\varepsilon}{2}.
\end{align*}
For fixed $j$, there exists $N\in\mathbb{N}$ such that
\begin{align*}
\underset{t\in[0,T]}{\sup}
|\langle h^n(t)-h^{\infty}(t),S_j\phi\rangle|
<\frac{\varepsilon}{2},
\qquad \forall\, n>N.
\end{align*}
Hence, we conclude that
$h^n\to h^{\infty}$ in
$C_w([0,T];B^s_{p,\infty})$.

In conclusion, combining with \textbf{Step 1--Step 4}, we complete the proof of Theorem \ref{LW}.
\end{proof}

\subsection {The proof of wave breaking}\label{sec4}
In this section, without relying on any conservation law, we establish an a priori estimate for the $L^\infty$-norm of strong solutions to Eq.~\eqref{dxfc} and use it to prove the wave-breaking criterion stated in Theorem~\ref{th.4}.
\begin{lemma}\label{Non es}
Suppose $h_0\in H^s(\mathbb{R})$ with $s>3/2$, and the maximal existence time $T^*>0$ be the lifespan of Cauchy problem \eqref{cauchy}. Then we have
\begin{align}\label{xygj}
    \|h(t)\|_{L^\infty} \leq \frac{1}{2 - \frac{4\|h_0\|_{L^\infty}}{2\|h_0\|_{L^\infty} + 1} \mathrm{e}^{\frac{t}{2}}} - \frac{1}{2},  \quad t \in \left[ 0, 2 \ln \left( 1 + \frac{1}{2\|h_0\|_{L^\infty}} \right) \right) \cap [0, T^*).
\end{align}
\end{lemma}

\begin{proof}
The characteristic flow associated with Eq.~\eqref{dxfc} is defined by
\begin{align}\label{tz}
\begin{cases}
    q_t(t,x)=\big(\frac{3}{2}h-\frac{1}{2}\big)(t,q(t,x)), & \hbox{$t\in[0,T^*),x\in\mathbb{R}$;} \\
    q(0,x)=x, & \hbox{$x\in\mathbb{R}.$}
\end{cases}
\end{align}
By the classical theory of ordinary differential equations, there exists a
unique solution $q\in C^1([0,T^*)\times\mathbb{R};\mathbb{R})$ to
\eqref{tz}. Moreover, for each $t\in[0,T^*)$, the map
$q(t,\cdot)$ is an increasing diffeomorphism of $\mathbb{R}$, and
\begin{align*}
q_x(t,x)=\exp{\left(\int_0^t\frac{3}{2}h_x(s,q(s,x))ds\right)}>0,\quad\forall (t,x)\in[0,T^*)\times\mathbb{R}.
\end{align*}
Consequently, the $L^\infty$ norm of any function $h(t,\cdot)\in L^{\infty}$ is preserved under $q$, that is
\begin{align*}
\|h(t,x)\|_{L^\infty}=\|h(t,q(t,x))\|_{L^\infty}.
\end{align*}
In view of \eqref{tz}, from Eq.~\eqref{dxfc}, we have
\begin{align}\label{eqcs}
\frac{d}{dt}[h(t,q(t,x))]&=\left(h_{t}+\big(\frac{3}{2}h-\frac{1}{2}\big)h_{x}\right)(t,q(t,x))\nonumber\\
&=\left(\frac{1}{2}Lh_x-\frac{1}{2}[L, Lh_x]h\right)(t,q(t,x)).
\end{align}
By using Young's inequality, we have
\begin{align}\label{ns1}
 \frac{1}{2}\left|Lh_x\right|=\frac{1}{2}|G_x*h|\leq\frac{1}{2}\|G_x\|_{L^1}\|h\|_{L^\infty}\leq \frac{1}{2}\|h\|_{L^\infty}   
\end{align}
and
\begin{align}\label{ns2}
\frac{1}{2}|[L, Lh_x]h|&\leq \frac{1}{2}|L((Lh_x)h)-(Lh_x)(Lh)|\nonumber\\
&\leq \frac{1}{2}|G*((Lh_x)h)|+\frac{1}{2}|(Lh_x)(Lh)|\nonumber\\
&\leq \frac{1}{2}\|G\|_{L^1}\|Lh_x\|_{L^\infty}\|h\|_{L^\infty}+\frac{1}{2}\|Lh_x\|_{L^\infty}\|Lh\|_{L^\infty}\nonumber\\
&\leq \frac{1}{2}\|G_x\|_{L^1}\|h\|_{L^\infty}^2+\frac{1}{2}\|G_x\|_{L^1}\|G\|_{L^1}\|h\|_{L^\infty}^2\nonumber\\
&\leq \frac{1}{2}\|h\|_{L^\infty}^2+\frac{1}{2}\|h\|_{L^\infty}^2\nonumber\\
&=\|h\|_{L^\infty}^2.
\end{align}
Combining \eqref{ns1} and \eqref{ns2},  from \eqref{eqcs}, it's easy to find that
\begin{align*}
    \left| \frac{dh(t, q(t,x))}{dt} \right| \leq  \|h(t)\|_{L^\infty}^2+\frac{1}{2}\|h(t)\|_{L^\infty}.
\end{align*}
Thus,
\begin{align*}
    - \|h(t)\|_{L^\infty}^2 -\frac{1}{2}\|h(t)\|_{L^\infty}\leq \frac{dh(t, q(t,x))}{dt} \leq  \|h(t)\|_{L^\infty}^2+\frac{1}{2}\|h(t)\|_{L^\infty}.
\end{align*}
Integrating the above inequality with respect to $t$ yields
\begin{gather*}
    - \int_0^t \left( \|h(\tau)\|_{L^\infty}^2 + \frac{1}{2} \|h(\tau)\|_{L^\infty} \right) d\tau + h_0(x) \leq h(t, q(t,x)), \\
    h(t, q(t,x)) \leq \int_0^t \left( \|h(\tau)\|_{L^\infty}^2 + \frac{1}{2} \|h(\tau)\|_{L^\infty} \right) d\tau + h_0(x).
\end{gather*}
Since $q(t,x)$ is a diffeomorphism of $\mathbb{R}$, by the above inequality, we obtain
\begin{align*}
    \|h(t)\|_{L^\infty} \leq \int_0^t \left( \|h(\tau)\|_{L^\infty}^2 + \frac{1}{2} \|h(\tau)\|_{L^\infty} \right) d\tau + \|h_0\|_{L^\infty}.
\end{align*}
Then by using Lemma \ref{Osgood}, we have
\begin{align*}
    \ln \frac{\|h(t)\|_{L^\infty}}{\|h(t)\|_{L^\infty} + \frac{1}{2}} \leq \ln \frac{\|h_0\|_{L^\infty}}{\|h_0\|_{L^\infty} + \frac{1}{2}} + \frac{1}{2}t,
\end{align*}
which implies
\begin{align*}
    \|h(t)\|_{L^\infty} \leq \frac{\frac{1}{2}}{1 - \frac{\|h_0\|_{L^\infty}}{\|h_0\|_{L^\infty} + \frac{1}{2}} \mathrm{e}^{\frac{1}{2}t}} - \frac{1}{2}.
\end{align*}
Indeed,
\begin{align*}
    t < 2 \ln \left( 1 + \frac{1}{2\|h_0\|_{L^\infty}} \right).
\end{align*}
Consequently, this implies the desired a priori estimate \eqref{xygj}, which completes the proof of Lemma \ref{Non es}.
\end{proof}

With Lemma \ref{Non es} in hand, we are now in a position to give the proof of Theorem~\ref{th.4}.

\begin{proof}[\normalfont\bfseries Proof of Theorem~\ref{th.4}]
Differentiating Eq.~\eqref{dxfc} with respect to $x$, we have
\begin{align}\label{eqxx}
h_{tx}+\big(\frac{3}{2}h-\frac{1}{2}\big)h_{xx}=-\frac{3}{2}h_x^2+\frac{1}{2}(G*h-h)-\frac{1}{2}\partial_x[L,Lh_x]h.
\end{align}
Notice that
\begin{align}\label{jhzgj}
\frac{1}{2}|\partial_x[L,Lh_x]h|
&\leq \frac{1}{2}\|\partial_xL((Lh_x)h)-\partial_x((Lh_x)(Lh))\|_{L^\infty}\nonumber\\
&\leq \frac{1}{2}\left(\|G_x*((Lh_x)h)\|_{L^\infty}+\|\partial_x^2(Lh)\|_{L^\infty}\|Lh\|_{L^\infty}+\|Lh_x\|_{L^\infty}^2\right)\nonumber\\
&\leq \frac{1}{2}\left(\|(Lh_x)h\|_{L^\infty}+(\|Lh\|_{L^\infty}+\|h\|_{L^\infty})\|Lh\|_{L^\infty}+\|Lh_x\|_{L^\infty}^2\right)\nonumber\\
&\leq \frac{1}{2} \left(\|h\|_{L^\infty}^2+2\|h\|_{L^\infty}^2 +\|h\|_{L^\infty}^2     \right)\nonumber\\
&\leq 2\|h\|_{L^\infty}^2,
\end{align}
where we have used Young's inequality and $\partial_x^2 (Lf) = Lf - f$.

Then, along the flow $q(t,x)$, combining \eqref{eqxx} and \eqref{jhzgj}, we obtain
\begin{align}\label{eq2}
\frac{d}{dt}[h_x(t,q(t,x_0))]&=\left[h_{tx}+\big(\frac{3}{2}h-\frac{1}{2}\big)h_{xx}\right](t,q(t,x_0))\nonumber\\
&=\left[-\frac{3}{2}h_x^2+\frac{1}{2}(G*h-h)-\frac{1}{2}\partial_x[L, Lh_x]h\right](t,q(t,x_0))\nonumber\\
&\leq -\frac{3}{2}h_x^2(t,q(t,x_0))+\frac{1}{2}\|G*h\|_{L^\infty}+\frac{1}{2}\|h\|_{L^\infty}+2\|h\|_{L^\infty}^2\nonumber\\
&\leq -\frac{3}{2}h_x^2(t,q(t,x_0))+2\|h\|_{L^\infty}^2+\|h\|_{L^\infty}\nonumber\\
&= -\frac{3}{2}h_x^2(t,q(t,x_0))+2\left(\|h\|_{L^\infty}+\frac{1}{4}    \right)^2-\frac{1}{8}\nonumber\\
&\leq -\frac{3}{2}h_x^2(t,q(t,x_0))+2\left(\frac{1}{2 - \frac{4\|h_0\|_{L^\infty}}{2\|h_0\|_{L^\infty} + 1} \mathrm{e}^{\frac{T_0}{2}}} - \frac{1}{4}  \right)^2-\frac{1}{8}.
\end{align}
Define
\begin{align*}
&m(t):=h_x(t,q(t,x_0)),\quad K:=2\left(\frac{1}{2 - \frac{4\|h_0\|_{L^\infty}}{2\|h_0\|_{L^\infty} + 1} \mathrm{e}^{\frac{T_0}{2}}} - \frac{1}{4}  \right)^2-\frac{1}{8}.
\end{align*}
Then, by \eqref{eq2}, for all $t\in [0,T_0]\cap[0,T^*)$, we have
\begin{align}
\frac{d m(t)}{dt}\leq -\frac{3}{2}m^2(t)+K.
\end{align}
The initial assumption \eqref{czjs} implies that
\begin{align*}
m(0)
&\leq -\frac{2\sqrt{\frac{2K}{3}}}{\mathrm{e}^{\sqrt{6K}\,T_0}-1}-\sqrt{\frac{2K}{3}}
<-\sqrt{\frac{2K}{3}}.
\end{align*}
Since
\begin{align*}
m'(t)
&\le -\frac32 m^2(t)+K \\
&=-\frac32\Bigl(m(t)-\sqrt{\tfrac{2K}{3}}\Bigr)\Bigl(m(t)+\sqrt{\tfrac{2K}{3}}\Bigr),
\end{align*}
by a continuity argument, we obtain
\begin{align*}
m(t)\leq m(0)<-\sqrt{\frac{2K}{3}},
\qquad t\in[0,T_0]\cap[0,T^*).
\end{align*}
Hence
\begin{align*}
-\left(\frac{m(t)+\sqrt{\frac{2K}{3}}}{m(t)-\sqrt{\frac{2K}{3}}}\right)'
\leq-\sqrt{6K}\,\frac{m(t)+\sqrt{\frac{2K}{3}}}{m(t)-\sqrt{\frac{2K}{3}}}.
\end{align*}
Integrating the above inequality over \([0,t]\), we get
\begin{align*}
\frac{m(0)+\sqrt{\frac{2K}{3}}}{m(0)-\sqrt{\frac{2K}{3}}}
\le
\frac{m(t)+\sqrt{\frac{2K}{3}}}{m(t)-\sqrt{\frac{2K}{3}}}
\mathrm{e}^{-\sqrt{6K}\,t}.
\end{align*}
Equivalently,
\begin{align*}
\frac{m(0)+\sqrt{\frac{2K}{3}}}{m(0)-\sqrt{\frac{2K}{3}}}
\mathrm{e}^{\sqrt{6K}\,t}-1
\leq\frac{2\sqrt{\frac{2K}{3}}}{m(t)-\sqrt{\frac{2K}{3}}}\leq 0.
\end{align*}
Since
\begin{align*}
\frac{m(0)-\sqrt{\frac{2K}{3}}}{m(0)+\sqrt{\frac{2K}{3}}}
>1,
\end{align*}
Therefore, there exists some $T^*$ satisfying
\begin{align*}
T^*&\leq\frac{1}{\sqrt{6K}}\ln\left(\frac{m(0)-\sqrt{\frac{2K}{3}}}{m(0)+\sqrt{\frac{2K}{3}}}\right)\leq T_0
\end{align*}
such that
\begin{align*}
\underset{t\rightarrow T^*}{\lim}m(t)=-\infty.
\end{align*}
This completes the proof of Theorem \ref{th.4}.
\end{proof}

\subsection{The proof of ill-posedness}\label{sec5}
In this section, motivated by the construction in \cite{r202,r505}, we prove that the Cauchy problem \eqref{cauchy} is ill-posed in
$B^{1+\frac{1}{p}}_{p,r}$ with $1\leq p\leq\infty$ and $1<r\leq\infty$ in the sense of norm inflation. To this end, we first establish an additional blow-up criterion for the Cauchy problem \eqref{cauchy} in the critical space.

\begin{lemma}\label{crbp}
Let $h_0\in H^s(\mathbb{R})$ with $s>\frac{3}{2}$, and let $T^*$ be the maximal existence time of the corresponding solution $h$ to Cauchy problem \eqref{cauchy}, then $h(t,x)$ blows up in finite time $T^*<\infty$ if and only if
\begin{align}\label{crcr}
\int_0^{T^*}
\|h(t')\|_{B_{\infty,\infty}^{1}}\,\mathrm{d}t'=\infty.
\end{align}
\end{lemma}

\begin{proof}
Applying $\Lambda^{s}$ to \eqref{cauchy} and taking the $L^2$ inner product with $\Lambda^{s}h$ yield 
\begin{align*}
\frac{1}{2}\frac{\mathrm{d}}{\mathrm{d}t}\|h\|_{H^s}^{2}+\frac{3}{2}
\left\langle\Lambda^s(hh_x),\Lambda^sh
\right\rangle={}&\frac{1}{2}\left\langle
\Lambda^sLh_x,\Lambda^sh\right\rangle+\frac{1}{2}\left\langle\Lambda^sh_x,\Lambda^sh\right\rangle-\frac{1}{2}
\left\langle\Lambda^s\big([L,Lh_x]h\big),\Lambda^sh\right\rangle.
\end{align*}
By Lemma~\ref{KP} and Cauchy-Schwarz inequality, we have
\begin{align*} \left| \left\langle\Lambda^s(hh_x),\Lambda^sh\right\rangle \right| &\leq \left\|[\Lambda^s,h]h_x\right\|_{L^2} \|h\|_{H^s} +\frac{1}{2}\|h_x\|_{L^\infty}\|h\|_{H^s}^2\\ &\leq C\|h_x\|_{L^\infty}\|h\|_{H^s}^2 \end{align*}
and
\begin{align*}
\left|\left\langle
\Lambda^s\big([L,Lh_x]h\big),\Lambda^sh\right\rangle\right|
&\leq\left\|\Lambda^s\big([L,Lh_x]h\big)
\right\|_{L^2}\|\Lambda^sh\|_{L^2}\\
&=\|[L,Lh_x]h\|_{H^s}\|h\|_{H^s}.
\end{align*}
Noticing that $\left\langle \Lambda^sLh_x,\Lambda^sh\right\rangle=
\left\langle\Lambda^sh_x,\Lambda^sh\right\rangle
=0$, it follows that
\begin{align}\label{cr1}
\|h(t)\|_{H^s}\leq{}&\|h_0\|_{H^s}
+C\int_0^t\left(
\|h_x(t')\|_{L^\infty}\|h(t')\|_{H^s}+\|[L,Lh_x]h(t')\|_{H^s}\right)\,\mathrm{d}t'.
\end{align}
By Lemma \ref{Besov}, Lemma \ref{product}, we have
\begin{align}\label{cr2}
\|[L,Lh_x]h\|_{H^s}
&\leq \|L((Lh_x)h)\|_{H^s}+\|(Lh_x)Lh\|_{H^s}\nonumber\\
&\leq C\big(\|Lh_x\|_{L^\infty}\|h\|_{H^s}
+\|h\|_{L^\infty}\|Lh_x\|_{H^s}+\|Lh_x\|_{L^\infty}\|Lh\|_{H^s}
+\|Lh\|_{L^\infty}\|Lh_x\|_{H^s}\big)\nonumber\\
&\leq C\|h\|_{L^\infty}\|h\|_{H^s}.
\end{align}
Combining \eqref{cr1} and \eqref{cr2}, we obtain
\begin{align}\label{cr3}
\|h(t)\|_{H^s}
\leq\|h_0\|_{H^s}
+C\int_0^t\left(
\|h_x(t')\|_{L^\infty}
+\|h(t')\|_{L^\infty}
\right)
\|h(t')\|_{H^s}\,\mathrm{d}t'.
\end{align}
By Lemma \ref{Besov} and the embedding $B_{\infty,1}^{1}\hookrightarrow W^{1,\infty}$, we have
\begin{align*}
\|h_x\|_{L^\infty}+\|h\|_{L^\infty}
\leq C\left[1+\|h\|_{B_{\infty,\infty}^{1}}
\log\big(\mathrm{e}+\|h\|_{H^s}\big)
\right].
\end{align*}
Therefore, \eqref{cr3} yields
\begin{align}
\|h(t)\|_{H^s}\leq{}&\|h_0\|_{H^s}
+C\int_0^t\left[1+\|h(t')\|_{B_{\infty,\infty}^{1}}\log\big(\mathrm{e}+\|h(t')\|_{H^s}\big)
\right]\|h(t')\|_{H^s}
\,\mathrm{d}t'.
\end{align}
Applying Gronwall's inequality, we obtain
\begin{align*}
\|h(t)\|_{H^s}
\leq\|h_0\|_{H^s}\exp\left\{
Ct+C\int_0^t
\|h(t')\|_{B_{\infty,\infty}^{1}}
\log\big(\mathrm{e}+\|h(t')\|_{H^s}
\big)\,\mathrm{d}t'
\right\}.
\end{align*}
Taking the logarithm of the above inequality and applying Gronwall's inequality again yields
\begin{align}\label{cr4}
\log\big(\mathrm{e}+\|h(t)\|_{H^s}
\big)\leq\left[\log\big(\mathrm{e}+\|h_0\|_{H^s}
\big)
+Ct\right]\exp\left\{
C\int_0^t\|h(t')\|_{B_{\infty,\infty}^{1}}
\,\mathrm{d}t'
\right\}.
\end{align}
If $T^*$ is finite and $\int_0^{T^*}\|h(t')\|_{B_{\infty,\infty}^{1}}\,\mathrm{d}t'<\infty$, then \eqref{cr4} implies that $\underset{{0\leq t<T^*}}{\sup}\|h(t)\|_{H^s}<\infty$, which contradicts the assumption that $T^*$ is the maximal existence time.

On the other hand, if $\int_0^{T^*}\|h(t')\|_{B_{\infty,\infty}^{1}}\,\mathrm{d}t'=\infty$ and the fact that $W^{1,\infty}\hookrightarrow B_{\infty,\infty}^{1}$, we conclude that
$h$ blows up in finite time. The proof is complete.
\end{proof}

\begin{proof}[\normalfont\bfseries Proof of Theorem~\ref{ill}]
Fix $1\leq p\leq \infty$ and $1<r\leq \infty$. Let $\widetilde{\varphi}\in\mathcal{D}(\mathbb{R})$ be an even, non-zero,
non-negative function such that $\widetilde{\varphi}\varphi=\widetilde{\varphi}$. For $N\in\mathbb{N}^+$ large enough, we define
\begin{align*}
h_0(x)=\sum_{n=2}^{N}\dfrac{\psi_n(x)}
{N^{1/r}(\ln N)2^{2n}},
\end{align*}
where
$\mathcal{F}(\psi_n)(\xi)=i2^{-n}\xi
\widetilde{\varphi}(2^{-n}\xi)$.
Since $\|\psi_n\|_{L^p}\leq C2^{n(1-\frac1p)}$ and
$\psi_{n,x}(0)=-c_0 2^{2n}$ for some $c_0>0$, it is easy to find that
$h_0\in H^\infty(\mathbb{R})\cap
B^{1+\frac{1}{p}}_{p,r}(\mathbb{R})$ is odd and satisfies
\begin{align*}
&\|h_0\|_{B_{p,r}^{1+\frac1p}}
\leq
\frac{C}{N^{1/r}\ln N}
\left\|\left(2^{n(1+\frac1p)}
\frac{\|\psi_n\|_{L^p}}{2^{2n}}\right)_{n=2}^{N}
\right\|_{\ell^r}
\leq\frac{C}{\ln N},\\
&\|h_0\|_{L^\infty}
\leq\frac{C}{N^{1/r}\ln N}
\sum_{n=2}^{N}2^{-n}
\leq\frac{C}{\ln N},\\
&h_{0,x}(0)
=-\frac{c_0(N-1)}{N^{1/r}\ln N}
\leq-\frac{CN^{1-\frac1r}}{\ln N}.
\end{align*}
Taking $x_0=0$ and $T_0=(\ln N)^{-1}$, for $N$ sufficiently large, one can verify that
\begin{align*}
0<T_0<2\ln\left(1+\frac{1}{2\|h_0\|_{L^\infty}}\right).
\end{align*}
Moreover, the constant $K$ defined in Theorem~\ref{th.4}
satisfies
\begin{align*}
0<K=2\left(\frac{1}{2-\dfrac{4\|h_0\|_{L^\infty}}
{2\|h_0\|_{L^\infty}+1}\mathrm{e}^{\frac{T_0}{2}}}-\frac{1}{4}\right)^2-\frac{1}{8}=
\frac{\dfrac{4\|h_0\|_{L^\infty}}
{2\|h_0\|_{L^\infty}+1}
\mathrm{e}^{\frac{T_0}{2}}}{\left(2-\dfrac{4\|h_0\|_{L^\infty}}
{2\|h_0\|_{L^\infty}+1}\mathrm{e}^{\frac{T_0}{2}}\right)^2}\leq\frac{C}{\ln N},
\end{align*}
and
\begin{align*}
h_{0,x}(0)\leq -\ln N\leq
-\frac{2\sqrt{\frac{2K}{3}}}{\mathrm{e}^{\sqrt{6K}\,T_0}-1}-\sqrt{\frac{2K}{3}}.
\end{align*}
Therefore, all the assumptions of Theorem~\ref{th.4} are
satisfied. Hence, the corresponding solution $h\in C([0,T);H^\infty(\mathbb{R}))$ blows up in finite time and $0<T\leq T_0=\frac{1}{\ln N}$. Then, by Lemma \ref{crbp}, we have
\begin{align*}
\int_0^T\|h(t)\|_{B^1_{\infty,\infty}}\,\mathrm{d}t=\infty.
\end{align*}
Therefore,
\begin{align*}
\|h\|_{L^\infty\left([0,T);B^{1+\frac{1}{p}}_{p,r}\right)}
\geq C\|h\|_{L^\infty\left([0,T);B^1_{\infty,\infty}\right)}
=\infty,
\end{align*}
where we used the embedding
$B_{p,r}^{1+\frac1p}
\hookrightarrow B_{\infty,r}^{1}
\hookrightarrow B_{\infty,\infty}^{1}$.
Thus, we obtain
\begin{align*}
\|h\|_{L^\infty\left([0,T);B^{1+\frac{1}{p}}_{p,r}\right)}
\geq\ln N.
\end{align*}
This completes the proof.
\end{proof}

\section*{Statements and Declarations}

\smallskip
\noindent\textbf{Funding}

\noindent
Yang was supported by the National Natural Science
Foundation of China (Grant No. 12561101).
Yin was supported by the National Natural Science
Foundation of China (Grant No. 12571261).

\smallskip
\noindent\textbf{Competing Interests}

\noindent
The authors declare that they have no conflict of interest.

\smallskip
\noindent\textbf{Data Availability}

\noindent
No data was used for the research described in the article.

\end{document}